\documentclass[11pt,a4paper,reqno]{article} 
\usepackage{amsthm,authblk} 
\newcommand{\address}[1]{\affil{#1}} 
\usepackage{fullpage}
\usepackage{amsmath,amssymb,amsfonts,latexsym,hyperref}
\usepackage{accents,graphicx}
\usepackage{blkarray}
\usepackage{color}
\usepackage{caption}
\usepackage{subcaption}
\newcommand{\sn}{S_n}
\newcommand{\sk}{S_k}

\newcommand{\ut}[1]{\underaccent{\tilde}{#1}}
\renewcommand{\vec}[1]{\ut{#1}}
\newcommand{\inv}{\mathrm{inv}}

\newcommand{\des}{\mathrm{des}}
\newcommand{\destat}{\mathrm{cddes}} % stands for CHATTERJEE-DIACONI
\newcommand{\cddes}{\destat}
\newcommand{\id}{\mathrm{id}}

\newcommand{\intenergy}{\phi}
\def\H{\mathcal{H}}

\newcommand{\stat}{\mathsf{stat}}

\newcommand{\Stat}{\mathsf{Stat}}

\newcommand{\smax}{s_{max}}
\newcommand{\CDdes}{\mathsf{CDdes}}

\newcommand{\GETOUT}[1]{}

\newcommand*{\doi}[1]{\href{http://dx.doi.org/#1}{doi:#1}}
\renewcommand{\P}{\mathbb{P}}

\newcommand{\PSC}{\mathsf{PSC}}

\newtheorem{theorem}{Theorem} % [section]
\newtheorem{proposition}[theorem]{Proposition}

\newtheorem{comment}[theorem]{Comment}

\theoremstyle{definition}

\newtheorem*{remark*}{Remark}

\newtheorem*{example*}{Example}

\title{An Ising model having permutation spin motivated by a permutation complexity measure}
\author{Mark Dukes} \address{School of Mathematics and Statistics\\ University College Dublin\\ Dublin 4, Ireland} \date{\it Dedicated to Einar Steingr\'imsson on the occasion of his retirement.}

\begin{document}
\maketitle
\begin{abstract}
In this paper we define a variant of the Ising model in which spins are replaced with permutations. The energy between two spins is a function of the relative disorder of one spin, a permutation, to the other. This model is motivated by a complexity measure for declarative systems. For such systems a state is a permutation and the permutation sorting complexity measures the average sequential disorder of neighbouring states. To measure the relative disorder between two spins we use a symmetrized version of the descent permutation statistic that has appeared in the works of Chatterjee \& Diaconis and Petersen. The classical Ising model corresponds to the length-2 permutation case of this new model. We consider and prove some elementary properties for the 1D case of this model in which spins are length-3 permutations. 
\end{abstract}
\section{Introduction}
\label{sectionone}
Since it was introduced by Bandt and Pompe~\cite{bp} in 2002, permutation entropy has 
become an established time series tool~\cite{parlitz,li,zanin,ruiz,zuninoreview} and
been developed in several ways~\cite{amigo}.
The permutation entropy of a time series is formed by mapping the sequence to a sequence of contiguous length-$k$ blocks, 
and then transforming the entries in each of these blocks into a length-$k$ permutation that summarizes how the values in the blocks relate to one another. 
In the permutations literature this transformation is called the {\em standardization}. 
For example, the standardization of the length-4 block $(0.5,7.1,-1,2.75)$ is the permutation $(2,4,1,3)$.
To then calculate the permutation entropy of the time series one considers $k!$ empty boxes labelled with the permutations of the set $\{1,\ldots,k\}$.
For each permutation that appears as the standardization of a block, put a ball into the box having that permutation as its label.
The permutation entropy is then the Shannon entropy of the distribution of balls amongst boxes.

Central to the calculation of the permutation entropy is the representation of the data or sequence as a sequence of permutations:
$$(\pi^{(1)},\pi^{(2)},\ldots,\pi^{(m)}) \mbox{ where  each } \pi^{(i)} \in \sk,$$
and $\sk$ is the set of permutations of the set $\{1,\ldots,k\}$. 
Permutations representing states of a system is not a new concept~\cite{ayyer,hopkins,dsss,dsssferrers,einarlauren}.
Recently, permutations (of subsets) of a set have been used to classify the allowed outcomes of declarative processes~\cite{ijgs}.
Declarative processes are processes that can happen in any way so long as it is not forbidden by a set of declarative constraints.
A declarative constraint specifies some type of dependency on the order in which activities can occur.
The execution of a declarative process is a permutation of the activities showing the order in which they occur.
Such processes have been used to model policy and healthcare processes~\cite{dukesieee,dukesfdap}, and a comparative analysis of such processes has been achieved through the introduction of combinatorial diversity metrics for these systems.

When a permutation represents the state of a system, information regarding how different the subsequent state is to the current one 
is essential in order to gain an insight into the diversity of such a process.
In particular, if one considers the day-to-day execution of healthcare processes, then information that there are wildly differing executions or extremely consistent executions of such a process help to inform the planning needs and resource requirements of the practitioners and healthcare providers.

Although the permutation complexity mentioned at the beginning of this paper is not a measure that captures the properties in which we are currently interested, 
the intermediate representation of the time series as a sequence of permutations in the calculation of the permutation entropy is useful.
The simplest and most natural candidate for a permutation sorting complexity measure to capture the sequential disorder in relation to what we will call a {\it{permutation process}}:
 $$\Pi = (\pi^{(1)},\pi^{(2)},\ldots,\pi^{(m)}),$$ where each $\pi^{(i)} \in\sk$, is 
an averaging over all neighbouring states:
\begin{align}
\PSC^{(\textrm{lin})}(\Pi) ~=~ \PSC^{(\textrm{lin})}_{\intenergy,m}(\Pi) ~:=~ \dfrac{1}{m-1}\sum_{t=1}^{m-1} \intenergy(\pi^{(t)},\pi^{(t+1)}).
\end{align}
In this equation $\intenergy(\alpha,\beta)$ is a measure of how different the permutation $\beta$ is to $\alpha$. 
More will be mentioned about a precise form for $\intenergy$ in Section 2.
This complexity measure $\PSC^{(\textrm{lin})}$ can be seen to resemble a Hamiltonian for a yet-to-be-defined model on a line, with a yet-to-be specified interaction energy.
It is this fact that motivates the Ising model that we now define and study in this paper.

\subsection{A permaspin model}
The Ising model~\cite{ising25}, a well-established model of magnetism, has witnessed several variants~\cite{baxter} since its inception.
In this paper we introduce a toy model that is a variant of the Ising model that we term the {\em permaspin model}.
The key feature of this model that distinguishes it from other Ising-like models is that a site's spin is some permutation of the set $\{1,\ldots,k\}$.
The interaction energy between two neighbouring sites will be a measure of how unordered the state, a permutation, of one site is with respect to the state of the other site. 
The external field acting on a site will be a measure of the disorder of that site's spin, in and of itself.

First let us define the permaspin model on a general graph.
The length $k$ of a permutation representing the spin state of a site is a parameter of the model and we will signify this by using the term {\it{$k$-permaspin}}. 
Let $G=G(V,E)$ be a simple graph with vertex set $V=\{v_1,\ldots,v_n\}$ and edge set $E$.
Let $\sk$ be the set of all permutations of the set $\{1,\ldots,k\}$ and let $\id$ be the identity permutation.
It will be useful to use the notation $\sk^n = \sk \times \cdots \times\sk$ for the $n$-fold Cartesian product that will be the state space of the system.

Given a configuration of permaspins $\vec{\pi} := (\pi^{(1)},\ldots, \pi^{(n)})\in \sk^n$ on $G$ where $\pi^{(i)}$ is the permaspin of vertex $v_i$, 
let us define the energy of this configuration through the Hamiltonian
\begin{align}\label{general_hamiltonian}
\H(\vec{\pi}) = -J \sum_{(v_i,v_j) \in E} \intenergy(\pi^{(i)}, \pi^{(j)}) - H \sum_{v_i \in V} \intenergy^{\textrm{external}}(\pi^{(i)}),
\end{align}
for some yet to be specified interaction energies $\intenergy$ and $\intenergy^{\textrm{external}}$.
We consider $J$ to be positive so that the model is ferromagnetic.
Note that the Hamiltonian for the line graph with $J=-1$ and $H=0$ is the complexity measure $\PSC^{(\textrm{lin})}/|E|$.
The partition function for the permaspin model is
\begin{align} \label{partfuncteqn}
Z(\beta) =  \sum_{ \vec{\pi}  \in \sk^m } \exp (-\beta \H (\vec{\pi} ) ) ,
\end{align}
where $\H(\vec{\pi})$ is given in equation (\ref{general_hamiltonian}).
The free energy is 
\begin{align} \label{freeenergydef}
f(\beta) = -\lim_{n \to \infty} \dfrac{1}{\beta n} \log Z(\beta).
\end{align}
%XXXX
We will assume the configuration probability for this canonical ensemble is given by the Boltzmann distribution so that the probability that the state $\sigma$ has a prescribed sequence of permaspins $\vec{\pi}$ is
$\mathbb{P}_{\beta}(\sigma=\vec{\pi}) ~=~  {\exp (-\beta \H (\vec{\pi} ) )}/{Z(\beta)}.$

In this paper we will take some first steps in considering simple instances of the permaspin model.
In Section~\ref{sectiontwo} we will discuss the interaction energy for this model and make a choice for the interaction energy based on both interaction symmetry and simplicity.
Our choice of interaction energy yields the classical Ising model as the 2-permaspin model.
In Section~\ref{sectionthree} we consider the 1D $k$-permaspin model in the absence of an external field, and see how closed forms for both the partition function and free energy can be
given as simple transformations of the generating function for the interaction energy. 

In Section~\ref{sectionfour} we consider the 1D 3-permaspin model in an external field.
It is not possible to give closed form expressions for the six eigenvalues associated with this system.
However, we are able to rule out several eigenvalues from being largest and see that the characteristic polynomial for the associated transition matrix factors as a cubic and three linear terms (two of which are equal).
Experimentally, it seems the largest eigenvalue is one of the cubic roots and we note some properties in relation to it. 
We further note that one of the eigenvalues with an explicit expression could be a potentially close lower bound for the largest eigenvalue.
We consider a low-temperature approximation for the free energy of the 1D 3-permaspin model.
In Section~\ref{sectionsix} we consider the mean-field 1D 3-permaspin model that results in a sum for the partition function.
Finally, we conclude with a discussion in Section~\ref{sectionseven}.

\section{A double Eulerian interaction energy}
\label{sectiontwo}
For two neighbouring permaspins, $\sigma$ and $\pi$, we would like the interaction energy $\intenergy$ to be a measure of how much one differs from the other.
There are several ways to model this, but for the purposes of this paper we will consider this interaction energy to be a function of
the permutation $\tau$ that is required to sort $\sigma$ into $\pi$,
$$\tau = \sigma^{-1} \pi.$$
In the literature this permutation $\tau$ has become known as the {\it{transcript}} of the pair of permutations $(\sigma,\pi)$, see e.g.~\cite{transcript1,transcript2}.
If $\sigma$ and $\pi$ are the same then $\tau = \id$ whereas if $\sigma$ is the reverse of $\pi$ (i.e. $\sigma_i = \pi_{n+1-i}$) then $\tau$ will be $\id$ reversed. 
For permaspins $\alpha,\beta \in \sk$ and $\tau := \alpha^{-1} \beta$, suppose the generic permutation statistic $\stat(\tau)$ takes values in the set $\{0,1,\ldots,\smax\}$.
To normalize the interaction energy we define
$$\intenergy(\alpha,\beta) := 1-  \frac{2}{\smax}\stat(\alpha^{-1}\beta),$$
that takes values in the closed interval $[-1,1]\subset \mathbb{R}$.
A further requirement is one of symmetry: the interaction energy
between $\alpha$ and $\beta$ must be the same as the interaction energy between $\beta$ and $\alpha$, 
i.e. $\intenergy(\alpha,\beta) = \intenergy(\beta,\alpha)$. 
This is equivalent to the statistic $\stat$ satisfying 
\begin{equation}\label{statsym}
\stat(\pi) = \stat(\pi^{-1}) \mbox{ for all permutations $\pi$.}
\end{equation}

There are several permutation statistics that might be used to measure how unordered a given permutation is. 
These include, but are not limited to, the permutation statistics 
{\it number of descents} ($\des$), 
{\it number of inversions} ($\inv$), 
{\it number of exceedances}, and {\it{number of weak exceedances}}.
For example, the number of descents in a permutation $\pi=(\pi_1,\ldots,\pi_n)$ is the number of indices $i$ for which $\pi_i>\pi_{i+1}$.
A permutation statistic that is equidistributed with the descent statistic on permutations is called {\it{Eulerian}} whereas one that is equidistributed with the major index permutation statistic is called {\it{Mahonian}}.
Many of the more common permutation statistics are either {\it{Eulerian}} or {\it{Mahonian}}.
Different statistics will of course lead to different calculation challenges, with those for Mahonian statistics being more involved due mainly to the range of values the statistic may take.

The inversions statistic, being the only one mentioned above to satisfy the interaction symmetry property for permaspins (Equation~\ref{statsym}) would seem like the natural choice.
However, even for the 3-permaspin case, the calculation and comparative analysis of eigenvalues for the inv statistic case (to be discussed in the last section) appears difficult. 
The reason for this is the larger range of values that the inversion statistic can take.
With this point in mind, in this paper we choose the symmetrized statistic $\destat$ whereby:
\begin{equation}\label{destateqn}
\destat (\pi) ~:=~ \des(\pi) + \des(\pi^{-1}).
\end{equation}
Note that $\smax$ for $\destat$ on $S_n$ is $2(n-1)$.
This statistic was investigated in Chatterjee \& Diaconis~\cite{persi} and Petersen~\cite{petersen13}.
The coefficients of the generating function for this statistic have become known as the {\it{double Eulerian numbers}}~\cite[A298248]{oeis}.
See Table~\ref{smalldestat} for the first few double Eulerian numbers.

\begin{table}
$$\begin{array}{@{\quad}l@{\quad}|@{\quad}l@{\quad}} \hline
n & \multicolumn{1}{c}{\CDdes_n(x) \phantom{\displaystyle\int_0^1}} \\ \hline
1& 1 \\
2& 1  +x^2   \\
3& 1+4x^2 + x^4 \\
4& 1+10x^2+2x^3+10x^4 + x^6 \\
5& 1+20x^2+12x^3+54x^4+12x^5+20x^6+x^8 \\
6& 1+35x^2+42x^3+ 212x^4 + 140x^5+212x^6 +42x^7+35x^8 +x^{10} \\ \hline
\end{array}$$
\caption{The generating functions $\CDdes_n(x) = \sum_{\pi \in \sn} x^{\destat(\pi)}$ for small $n$.\label{smalldestat}}
\end{table}

The generating function for this statistic is (see \cite[Theorem 2]{petersen13})
\begin{equation}\label{CDdesgf}
\mathsf{CDdes}_n(u) = \sum_{ \pi \in \sn} u^{\destat(\pi)} = (1-u)^{2n+2} \sum_{i,j\geq 1} \binom{ij+n-1}{n} u^{i+j-2}.
\end{equation}

We will write $\intenergy_{\destat}$ to signify the use of Equation~\ref{destateqn} in the interaction energy.
With the interaction energy defined using the $\destat$ statistic, the 2-permaspin model is the classical Ising model where the spins take values in $\{-1,+1\}$. 
This is easily seen via the identification of the permutations $(1,2)$ and $(2,1)$ with the classical spins $+1$, $-1$, respectively.
If $\alpha,\beta \in S_2$ with corresponding classical spins $s_{\alpha},s_{\beta}$, then one finds $\intenergy(\alpha,\beta) = s_{\alpha} s_{\beta}$.

One outstanding matter is the term in the Hamiltonian that represents the energy interaction $\intenergy^{\textrm{external}}$ of a permaspin with an hypothesized external field.
The term in the Hamiltonian (Equation~\ref{general_hamiltonian}) that corresponds to this is $-H \intenergy^{\textrm{external}}(\pi)$. %, even though $\intenergy(\pi)$ as defined above is a function of two permutations.
We will assume that $\intenergy^{\textrm{external}}(\pi) := 1-  \frac{2}{\smax}\stat(\pi)$, which is precisely the same value one gets from assuming this is the energy of the permaspin $\pi$ with another permaspin that is the identity permutation, i.e. $\intenergy^{\textrm{external}}(\pi) = \intenergy(\id,\pi)=\intenergy(\pi,\id)$. 
In other words, one can consider the hypothesized external field to have permaspin $\id$.

%%%%%%%%%%%%%%%%%%%%%%%%%%%%%%%%% %%%%%%%%%%%%%%%%%%%%%%%%%%%%%%%%% %%%%%%%%%%%%%%%%%%%%%%%%%%%%%%%%% %%%%%%%%%%%%%%%%%%%%%%%%%%%%%%%%%
\section{The 1D permaspin model without external field} \label{sectionthree}
In this section we consider the 1D permaspin model without external field.
There are two distinct cases to consider: the path graph and the cycle graph.
The path graph is the graph with vertex set $V(G)=\{v_1,\ldots,v_n\}$ and edge set $E(G) = \{ (v_i,v_{i+1}) ~:~ 1 \leq i <n\}$. %, and the convention $v_{n+1}:=v_1$.
The cycle graph is the path graph with vertex $v_n$ joined to $v_1$.
The partition function for this 1D model is easily calculated since the partition function can be written as a deformation of the generating function for the $\stat$ statistic.

\begin{proposition}\label{1dnoexternalstat}
The partition functions for the path and cycle graphs are:
\begin{align*}
Z^{(\stat,path)}_n(\beta) = & k! e^{\beta J} \left( \exp(\beta J)  \Stat_k\left( 
									\exp\left( \frac{-2\beta J}{\smax} \right) 
									\right)\right)^{n-1}, \mbox{ and}\\
Z^{(\stat,cycle)}_n(\beta) = & \left( \exp(\beta J)  \Stat_k\left( 
									\exp\left( \frac{-2\beta J}{\smax} \right) 
									\right)\right)^{n},
\end{align*}
where $\Stat_k(x) := \sum_{\pi \in \sk} x^{\stat(\pi)}$.
The free energy for both cases is
$$f^{(\stat)}(\beta) = -J 
	-\dfrac{1}{\beta} \ln  
							\Stat_k\left(
								  e^{-2\beta J/\smax} 
                			\right).
$$
\end{proposition}

\begin{proof}
As the the edges are pairs $(v_i,v_{i+1})$ the partition function is 
\begin{align*}
Z^{(\stat,path)}_n(\beta)&= \sum_{ \vec{\pi}  \in \sk^n } \exp \left(\beta J \sum_{i=1}^{n-1}  \intenergy(\pi^{(i)}, \pi^{(i+1)})  \right) \\
&= \sum_{ \vec{\pi}  \in \sk^n } \exp \left(\beta J \sum_{i=1}^{n-1}  
	\left(
	1-\dfrac{2}{\smax} \stat((\pi^{(i)})^{-1} \pi^{(i+1)} )
	\right)
  	\right).
\end{align*}
Unity from the internal sum can be extracted and moved outside, and the internal exponent containing a sum can be written as a product of exponents:
\begin{align*}
Z^{(\stat,path)}_n(\beta) &= \exp( \beta J n) \sum_{ \vec{\pi}  \in \sk^n } \exp \left(\dfrac{-2\beta J}{\smax} \sum_{i=1}^{n-1}  
	\left(
	\stat((\pi^{(i)})^{-1} \pi^{(i+1)} )
	\right)
  	\right) \\
&= \exp( \beta J n) \sum_{ \vec{\pi}  \in \sk^n } 
		\prod_{i=1}^{n-1} 
		\exp \left(\dfrac{-2\beta J}{\smax} 
		\stat((\pi^{(i)})^{-1} \pi^{(i+1)} )
  		\right). 
\end{align*}
Fix $\pi^{(1)}$ and set $\sigma^{(1)} = (\pi^{(1)})^{-1} \pi^{(2)}$. 
Instead of summing over all $\pi^{(2)} \in \sk$, we can exploit the form in the exponent to sum over all $\sigma^{(1)} \in \sk$, and extend this to 
$\sigma^{(2)},\ldots,\sigma^{(n-1)}$.
The above expression becomes:
\begin{align*}
Z^{(\stat,path)}_n(\beta) &= 
		\exp( \beta J n) 
		\sum_{\pi^{(1)} \in \sk}
		\sum_{\sigma^{(1)} \in \sk} e^{\frac{-2\beta J}{\smax} \stat(\sigma^{(1)})} 
		\sum_{\sigma^{(2)} \in \sk} e^{\frac{-2\beta J}{\smax} \stat(\sigma^{(2)})}
		\cdots
		\sum_{\sigma^{(n-1)} \in \sk} e^{\frac{-2\beta J}{\smax} \stat(\sigma^{(n-1)})}\\
	&=  \exp(\beta J n) k! \left( \Stat_k\left( \exp\left( \frac{-2\beta J}{\smax} \right) \right)\right)^{n-1}\\
	&=  k! e^{\beta J} \left( \exp(\beta J)  \Stat_k\left( 
									\exp\left( \frac{-2\beta J}{\smax} \right) 
									\right)\right)^{n-1},
\end{align*}
where $\Stat_k(x) := \sum_{\pi \in \sk} x^{\stat(\pi)}$.
For the cycle graph, the derivation is almost identical to the path case except, when writing the partition function as a sum over all $(\pi^{(1)},\sigma^{(1)},\ldots,\sigma^{(n)})$, 
the sum over all $\pi^{(1)} \in \sk$ does not appear but there is an extra term for the interaction between the edges $v_n$ and $v_{n+1}=v_1$ which was absent for the path graph case.
The partition function is
\begin{align*}
Z^{(\stat,cycle)}_n(\beta) &= 
		\exp( \beta J n) 
		\sum_{\sigma^{(1)} \in \sk} e^{\frac{-2\beta J}{\smax} \stat(\sigma^{(1)})} 
		\sum_{\sigma^{(2)} \in \sk} e^{\frac{-2\beta J}{\smax} \stat(\sigma^{(2)})}
		\cdots
		\sum_{\sigma^{(n)} \in \sk} e^{\frac{-2\beta J}{\smax} \stat(\sigma^{(n)})}\\
	&=  \exp(\beta J n)  \left( \Stat_k\left( \exp\left( \frac{-2\beta J}{\smax} \right) \right)\right)^{n}\\
	&=  \left( \exp(\beta J)  \Stat_k\left( 
									\exp\left( \frac{-2\beta J}{\smax} \right) 
									\right)\right)^{n}.
\end{align*}
The free energy for both cases is
\begin{align*}
f^{(\stat)}(\beta) = -\lim_n \frac{1}{\beta n} \ln Z^{(\stat)}_n(\beta) =  
	-\dfrac{1}{\beta} \ln \left( 
							\exp(\beta J)  
							\Stat_k\left(
								\exp\left( \frac{-2\beta J}{\smax} \right)
                			\right)
						\right).
\end{align*}
\end{proof}

For the case in which the length of permaspins is large, we have the following approximation for the free energy.

\begin{proposition}\label{proptwo}
Let $\stat$ be the double Eulerian statistic $\destat$ and consider the $k$-permaspin model for $k$ large. 
Then the free energy for both the path and cycle graphs is 
\begin{align*}
f^{(\destat)}(\beta) 
	\approx &  %-2J 
 -2J - \dfrac{k  \beta J^2}{12} - \dfrac{1}{\beta}(k \ln k -k + O(\ln k)).
\end{align*}
\end{proposition}

\begin{proof}
Chatterjee and Diaconis~\cite[Thm 1.1]{persi} established 
the following properties for 
the statistic $\destat$.
For $\pi$ chosen uniformly at random from the symmetric group $\sk$, 
$$\mathbb{E}(\destat(\pi))=k-1 \mbox{ and } \mathrm{var}(\destat(\pi))= \dfrac{k+7}{6} - \dfrac{1}{k},$$
for $k \geq 2$.
Normalized by its mean and variance, $\destat(\pi)$ has a limiting standard normal distribution.
When the permutation statistic $\stat$ is $\destat$, by Proposition~\ref{1dnoexternalstat} 
the free energy for both cases is
\begin{align*}
f^{(\destat)}(\beta) 
	= &  -J - \dfrac{1}{\beta} \ln  \mathsf{CDdes}_k\left( e^{-\beta J/(k-1)}  \right)
	=   -J - \dfrac{1}{\beta} \ln  \sum_{\pi \in \sk} e^{-\beta J \cddes(\pi)/(k-1)}  
\end{align*}
where $\mathsf{CDdes}_k(x)$ is given in Equation~\ref{CDdesgf}.
Let $z=e^{-\beta J}$ and consider
$$g_k(z) = \sum_{\pi \in \sk} z^{\cddes(\pi)/(k-1)}.$$
This is the probability generating function of the random variable $X_k := \frac{1}{k-1}\cddes_k$. 
As $\cddes_k$ has mean $k-1$ and variance $\dfrac{k+7}{6} - \dfrac{1}{k}$, the random variable $X_k$ will 
have mean $\mathbb{E}(X_k)=1$ and variance $\mathrm{var}(X_k) = \dfrac{k+7}{6(k-1)^2} - \dfrac{1}{k(k-1)^2}.$
Since $\cddes_k$, when  normalized by its mean and variance, has a limiting standard normal distribution, 
the same is true for $X_k$. 
The p.g.f. $$g_k(z) = k! \sum_{x \in \{0,1/(k-1),\ldots,2\}} \P(X_k=x) z^x$$
so for $k$ large:
\begin{align*}
g_k(e^{-\beta J}) 
\approx & k! \int  \dfrac{1}{\sigma \sqrt{2\pi}} \exp\left({-\frac12 \left( \frac{x-1}{\sigma}\right)^2}\right)  (\exp({-\beta J}))^x  dx \\
= & k! \int  \dfrac{1}{\sigma \sqrt{2\pi}} \exp\left({-\frac12 \left( \frac{x-(1+\sigma^2 \beta J)}{\sigma}\right)^2}\right)  \exp\left({-\frac{1}{2\sigma^2} (1-(1+\sigma^2\beta J)^2)}\right)  dx \\
= & k! \exp\left({\beta J + \frac12 \sigma^2\beta^2 J^2}\right).
\end{align*}
For $k$ large, the free energy is then
\begin{align*}
f^{(\destat)}(\beta) 
	\approx &  -J - \dfrac{1}{\beta} \ln \left(k! 
			\exp\left(\beta J + \frac12 \sigma^2\beta^2 J^2\right) 
			\right) \\
	= &  -J - \dfrac{1}{\beta} \left(k \ln k -k + O(\ln k) +  \beta J + \frac12 \sigma^2\beta^2 J^2 \right) \\
	= &  -2J - \dfrac{k  \beta J^2}{12} - \dfrac{1}{\beta}(k \ln k -k + O(\ln k)).\qquad\qquad\qedhere
\end{align*}
\end{proof}
Notice that the high and low temperature contributions in Proposition~\ref{proptwo} are clear:
$$f^{(\destat)}(\beta \mbox{ small}) \approx - \dfrac{1}{\beta}(k \ln k -k) \qquad \mbox{ and } \qquad f^{(\destat)}(\beta \mbox{ large})  \approx - \dfrac{k  \beta J^2}{12}.$$
As we mentioned in the introduction, the $\cddes_2$ interaction energy corresponds to the classical 1D Ising model and, using $\CDdes_2(x)=1+x^2$, one recovers the well-known free energy for this case:
$f^{(\destat)}(\beta) = -\frac{1}{\beta} \ln(2 \cosh ({\beta J})).  $

\section{The 1D 3-permaspin model in an external field}
\label{sectionfour}
In the previous section we saw that knowledge of the generating function for the $\cddes$ statistic and its limiting behaviour allowed us to present several exact and limiting results.
In this section we will consider the 1D 3-permaspin model in an external field. 
As the 2-permaspin model on a general graph corresponds to the classical Ising model on that same graph, any properties relating to the latter are also properties of the 2-permaspin model.

The 3-permaspin model is the simplest one to first consider and for which results are not known.
It is a well-known fact that (see e.g. \cite[Cor. 6.41]{friedli}) that classical lattice models with finitely-many spin states per vertex and local interactions cannot exhibit a phase transition.
In order to calculate the partition function and free energy for the 3-permaspin model on a ring, 
we require the eigenvalues of an order 6 matrix.
The partition function for the general 1D $3$-permaspin model is:
\begin{align*}
Z^{(\stat)}_n(\beta) 
&= \sum_{ \vec{\pi}  \in S_3^n } \exp \left(
		-\beta \left(  
			-J \sum_{(v_i,v_j)\mbox{\tiny n.n.} } \intenergy(\pi^{(i)}, \pi^{(j)})  
		 	-H \sum_{v_i}  \intenergy(\pi^{(i)})
			\right)\right) \\
&= \sum_{ \vec{\pi}  \in S_3^n } \exp \left(
					\beta J \sum_{i=1}^{n}  
						\left( 1-\dfrac{2 \stat((\pi^{(i)})^{-1} \pi^{(i+1)}) }{\smax} \right)
					+\beta H \sum_{i=1}^{n} 
						\left( 1-\dfrac{2 \stat(\pi^{(i)}) }{\smax} \right)
					\right) \\
& = e^{\beta (J+H)n} \sum_{ \vec{\pi}  \in S_3^n  } \exp \left( 
					-\dfrac{2\beta J}{\smax}  \sum_{i=1}^{n}  
						\left( \stat((\pi^{(i)})^{-1} \pi^{(i+1)})  \right)
					-\dfrac{2\beta H}{\smax}  \sum_{i=1}^{n} 
						\left( \stat(\pi^{(i)})  \right)
					\right).
\end{align*}
Set $J'= -2\beta J/\smax$ and $H'= -2\beta H /\smax$ so that the previous expression becomes:
\begin{align}
Z^{(\stat)}_n(\beta) 
&= e^{\beta (J+H)n} \sum_{ \vec{\pi}  \in S_3^n } \exp \left( 
						J' \sum_{i=1}^{n} \left( \stat((\pi^{(i)})^{-1} \pi^{(i+1)})  \right)
						+H' \sum_{i=1}^{n} \left( \stat(\pi^{(i)}  \right)
				\right)\nonumber \\
&= e^{\beta (J+H)n} \sum_{ \vec{\pi}  \in S_3^n  } 
					\prod_{i=1}^n  
						\exp\left(\dfrac{H' \stat(\pi^{(i)})}{2}
							+{J' \stat((\pi^{(i)})^{-1} \pi^{(i+1)}) } + 
						\dfrac{H' \stat(\pi^{(i+1)})}{2}\right) \nonumber \\
&= e^{\beta (J+H)n} \sum_{ \vec{\pi}  \in S_3^n  } 
					\prod_{i=1}^n  
						\exp\left(\dfrac{H' (\stat(\pi^{(i)}) + \stat(\pi^{(i+1)}))}{2}
							+{J' \stat((\pi^{(i)})^{-1} \pi^{(i+1)}) }\right) .
\end{align}

\begin{figure}
\centering
\begin{subfigure}{.5\textwidth}
  \centering
  \includegraphics[width=1.0\linewidth]{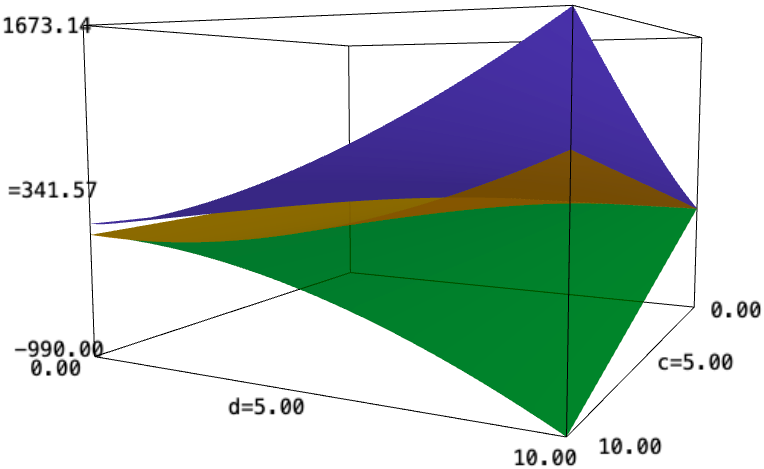}
  \caption{}
  \label{figblue1}
\end{subfigure}%
\begin{subfigure}{.5\textwidth}
  \centering
  \includegraphics[width=1.0\linewidth]{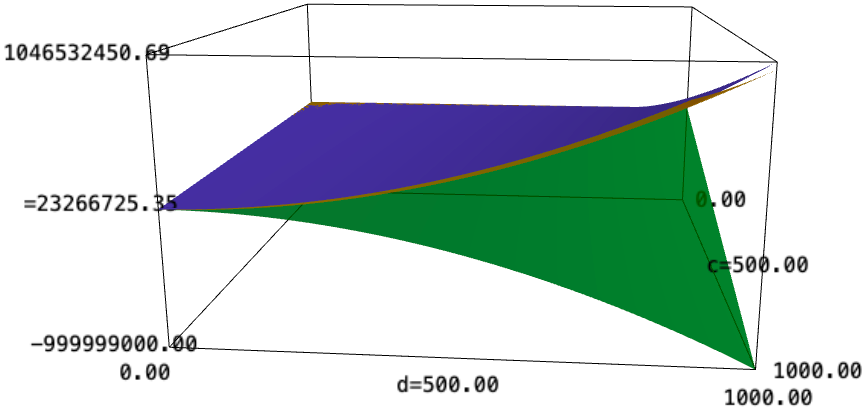}
\ \\[1.5em]
  \caption{}
  \label{figblue2}
\end{subfigure}
\caption{Surfaces for $\lambda^{*}$ (blue), $\lambda_4$ (orange),  and $\lambda_5$ (green).}
\label{threeeigens}
\end{figure}

Let $A$ be the $6 \times 6$ wherein row $i$th corresponds to the $i$th permutation of $S_3$ listed in lexicographic order, and the same labelling convention for the columns.
Define entry $A_{\pi,\sigma}$ of $A$ to be
\begin{align*}
A_{\pi,\sigma} := & 
\exp\left({\frac{ H'}{2} (\stat(\pi) + \stat(\sigma))  + J' \stat(\pi^{-1}\sigma)}\right)
= a^{\stat(\pi) + \stat(\sigma)} b^{\stat(\pi^{-1}\sigma)},
\end{align*}
where we have used the substitution $a= e^{\frac12 H'}$ and $b =e^{J'}$.
Using this matrix formalization, we can now write
\begin{align}\label{genpartitionfunction}
Z^{(\stat)}_n(\beta) &= e^{\beta (J+H)n} \mathrm{Tr}(A^n) = e^{\beta (J+H)n} (\lambda_1^n + \ldots+ \lambda_k^n),
\end{align}
where $\lambda_1,\ldots,\lambda_k$ are the eigenvalues of $A$.
Since we assume the statistic $\stat$ to satisfy the interaction symmetry condition (Equation~\ref{statsym}) the matrix $A$ will be symmetric. 
That $A$ is symmetric guarantees all eigenvalues of $A$ are real.
We have
\begin{align}
A=
\left(\begin{array}{rrrrrr}
1 & a^{2} b^{2} & a^{2} b^{2} & a^{2} b^{2} & a^{2} b^{2} & a^{4} b^{4} \\
a^{2} b^{2} & a^{4} & a^{4} b^{2} & a^{4} b^{4} & a^{4} b^{2} & a^{6} b^{2} \\
a^{2} b^{2} & a^{4} b^{2} & a^{4} & a^{4} b^{2} & a^{4} b^{4} & a^{6} b^{2} \\
a^{2} b^{2} & a^{4} b^{4} & a^{4} b^{2} & a^{4} & a^{4} b^{2} & a^{6} b^{2} \\
a^{2} b^{2} & a^{4} b^{2} & a^{4} b^{4} & a^{4} b^{2} & a^{4} & a^{6} b^{2} \\
a^{4} b^{4} & a^{6} b^{2} & a^{6} b^{2} & a^{6} b^{2} & a^{6} b^{2} & a^{8}
\end{array}\right).
\end{align}
The characteristic polynomial of $A$ is
\begin{align*}
c_A(\lambda) = &
-\bigg(
a^{12} ( b^{12} + 2 b^{10} - 7 b^{8} + 7 b^{4} - 2  b^{2} - 1) \\
& \qquad
+ a^{4} \lambda (1 - 3 a^{8} b^{4} - a^{4} b^{8}  + 2  a^{8} b^{2}  + a^{8}  + a^{4}  - 3 b^{4}  + 2  b^{2}  )\\
& \qquad
+\lambda^2( - a^{4} - a^{8} - a^{4} b^{4} - 2 \, a^{4} b^{2}  - 1)
+ \lambda^{3} 
\bigg)\\ & {\left(a^{4} b^{4} - 2 \, a^{4} b^{2} + a^{4} - \lambda\right)} {\left(a^{4} b^{4} - a^{4} + \lambda \right)}^{2}.
\end{align*}
Substitute $c$ for $a^4$ and $d$ for $b^2$ to write
\begin{align*}
c_A(\lambda) = &
-\bigg(
(d^2 + 4 d + 1) c^{3} (d + 1) (d - 1)^3
-\lambda c (cd^3 + 3c^2d + cd^2 + c^2 + cd + c + 3d + 1)(d - 1) \\
& \qquad
-\lambda^2 (cd^2 + c^2 + 2cd + c + 1)
+\lambda^3 \bigg) \left(cd^2-2cd+c-\lambda\right)  \left(cd^2-c+\lambda\right).
\end{align*}
This yields eigenvalues 
$\lambda_4 =  c(1-d)^2$, 
$\lambda_5 =  c(1-d^2)$
of $A$. 
Since $A$ has only real eigenvalues, the other three real eigenvalues $\lambda_1,\lambda_2,\lambda_3$ are solutions to $t(\lambda)=0$ where
\begin{align*}
t(\lambda) :=& \lambda^3 - (cd^2 + c^2 + 2cd + c + 1) \lambda^2   -c (cd^3 + 3c^2d + cd^2 + c^2 + cd + c + 3d + 1)(d - 1) \lambda \\ 
& +(d^2 + 4 d + 1) c^{3} (d + 1) (d - 1)^3\\
=& A'\lambda^3+B'\lambda^2+C'\lambda+D'.
\end{align*}
Let $\lambda^{*}$ be the largest root of $t(\lambda)=0$.
It remains to see which of the eigenvalues $\lambda^{*}$, $\lambda_4$, and $\lambda_5$ is largest.
Numerical investigations indicate that  $\lambda^{*}$ appears to be the largest and this is supported by the evidence presented by 
in Figures~\ref{figblue1} and \ref{figblue2}. (The blue surface corresponds to $\lambda^{*}$ while the orange and green surfaces correspond to $\lambda_4$ and $\lambda_5$, respectively.)
We summarize the results of these  observations and analysis as follows:
\begin{proposition}
\label{conjecturefive}
If $\lambda^{*} = \max(\lambda^{*},\lambda_4,\lambda_5)$ above, then 
the free energy for the 1D 3-permaspin model is 
\begin{align*}
f^{(\destat)}(\beta) 
= & -(J+H) - \dfrac{1}{\beta} \ln 
	\left(  {(c^{2} d + c^{2} + 2 \, c d + c + 1)} - 2 \mathbf{Re}(z)
\right)
\end{align*}
where $z$ is the unique root of the cubic
\begin{align*}
\lefteqn{\lambda^3 - (cd^2 + c^2 + 2cd + c + 1) \lambda^2}\\
&   -c (cd^3 + 3c^2d + cd^2 + c^2 + cd + c + 3d + 1)(d - 1) \lambda
  +(d^2 + 4 d + 1) c^{3} (d + 1) (d - 1)^3=0.
\end{align*}
for which $\mathbf{arg}(z) \in [2\pi/3,\pi]$ and where $c=e^{-2\beta H}$ and $d=e^{-2\beta J}$.
\end{proposition}
For larger values of $c$ and $d$ it seems that that $\lambda_4$ serves as a close lower bound to the $\lambda^{*}$.
This suggests that free energy for the 1D 3-permaspin model at low-temperature has behaviour:
\begin{align*}
f^{(\destat)}(\beta) \approx & -\lim_n \dfrac{1}{\beta n} \ln (e^{\beta(J+H)} \lambda_4^n) = - \dfrac{1}{\beta }(\beta(J+H)+ \ln \lambda_4) \\
=& -(J+H) - \dfrac{1}{\beta } \ln c(1-d)^2\\
=& -(J+H) - \dfrac{1}{\beta } \ln a^4(1-b^2)^2.
\end{align*}
As $a=e^{-\beta H/2}$ and $b=e^{-\beta J}$ this becomes
\begin{flalign*}
&& f^{(\destat)}(\beta) 
\approx & -(J+H) - \dfrac{1}{\beta }  \ln (e^{-\beta H/2})^4  (1-(e^{-\beta J})^2)^2 &&\\
&& =& -(J+H) - \dfrac{1}{\beta }  \ln e^{-2\beta H}  (1-e^{-2\beta J})^2 &&\\
&& =& -(J+H) + 2H - \dfrac{2}{\beta } \ln  (1-e^{-2\beta J}).
\end{flalign*}

\begin{comment}\label{commentseven}
The free energy for the 1D 3-permaspin model at low-temperature seems to have the behaviour:
\begin{align*}
f^{(\destat)}(\beta) 
\approx &
 -(J-H)  - \dfrac{2}{\beta } \ln  (1-e^{-2\beta J}) \approx -(J-H) - \dfrac{2}{\beta } \ln \beta .
\end{align*}
\end{comment}
Notice that as $\beta \to \infty$ in the above expression, $-(J-H) - \dfrac{2}{\beta } \ln \beta \to -(J-H)$.

Consider the configurations that represent the `largest' terms in the partition function for low temperatures, i.e. when $\beta$ is large.
These will be when all permaspins are the same. The contribution to the partition function $Z$ for configurations of this type is
\begin{align}
\sum_{\pi \in S_3} \exp \left( n\beta J \intenergy (\pi,\pi) + n \beta H \intenergy(\pi,\id)\right)
&= e^{n\beta J} \sum_{\pi \in S_3} \exp \left(n \beta H \intenergy(\pi,\id)\right) \nonumber\\
&= e^{n\beta J} \left( e^{-n\beta H} + 4 + e^{n\beta H} \right). \label{simpleone}
\end{align}
In the event that
$H=0$ the expression \ref{simpleone} is dominated by the term 
$ e^{n\beta J}$; 
$H>0$ 
it is dominated by the term
$e^{n\beta (J+H)}$; 
$H<0$ 
it is dominated by the term
$e^{n\beta (J-H)}$. 
These dominating terms imply the following low-temperature behavior:
\begin{comment}\label{commentsix}
A lower bound for free energy for the 1D 3-permaspin model at low-temperature seems to have the behaviour:
$$f^{(\cddes)}(\beta) ~\approx ~ -(J+|H|).$$
\end{comment}
This approximation coincides with the prediction of Comment~\ref{commentseven} for $H\leq 0$, but differs from it significantly for $H>0$. 

\section{A mean-field 3-permaspin model}
\label{sectionsix}

Let us suppose that for a 3-permaspin model, each of the $n$ permaspins has $q$ nearest-neighbors,
so that the 1D model corresponds to $q=2$ whereas the 2D model corresponds to $q=4$.
In this section we will consider a mean-field 3-permaspin model. 
For a mean-field model, the total field acting on site $v_i$ is 
$$H\intenergy(\pi^{(i)},\id) ~+~ \dfrac{qJ}{n-1} \sum_{j\neq i} \intenergy(\pi^{(i)},\pi^{(j)}).$$
Here, the sum over all the interactions of site $i$ with the $n-1$ sites that are not $i$ is averaged out by dividing the sum by $n-1$. The $q$ in the numerator reflects the multiplicity of this quantity with respect to the number of neighbours of site $i$.
This is equivalent to replacing the Hamiltonian with 
$$\H^{(mean)} (\vec{\pi}) = - \dfrac{qJ}{n-1} \sum_{
		i\neq j
		} \intenergy(\pi^{(i)},\pi^{(j)})  - H\sum_{i=1}^n \intenergy(\pi^{(i)},\id).$$
Note that the first sum is over all $\binom{n}{2}$ distinct pairs of vertices $\{v_i,v_j\}$.

Since every permaspin interacts with every other permaspin, we can simplify the above sum by considering the distribution of permaspins in $\vec{\pi}$.
Let us write $n_{\sigma}$ for the number of $\pi^{(i)}$ that equal $\sigma$. 
Then $$n_{123}+n_{132}+n_{213}+n_{231}+n_{312}+n_{321} = n.$$
Table~\ref{tablesix}(a) illustrates $\destat(\sigma^{-1}\pi)$ for all $\sigma$ and $\pi$, and the internal energy corresponding to a nearest neighbour pairing is given in Table~\ref{tablesix}(b).
\begin{table}[!h]
\begin{center}
\fbox{
\begin{tabular}{c@{\qquad}c}	
$\begin{pmatrix}
0 & 2 & 2 & 2 & 2 & 4 \\
2 & 0 & 2 & 4 & 2 & 2 \\
2 & 2 & 0 & 2 & 4 & 2 \\
2 & 4 & 2 & 0 & 2 & 2 \\
2 & 2 & 4 & 2 & 0 & 2 \\
4 & 2 & 2 & 2 & 2 & 0
\end{pmatrix}$
& 
$\begin{pmatrix}
1 & 0 & 0 & 0 & 0 & -1 \\
0 & 1 & 0 & -1 & 0 & 0 \\
0 & 0 & 1 & 0 & -1 & 0 \\
0 & -1 & 0 & 1 & 0 & 0 \\
0 & 0 & -1 & 0 & 1 & 0 \\
-1 & 0 & 0 & 0 & 0 & 1
\end{pmatrix}$ \\ & \\
(a) $\destat(\sigma^{-1}\pi)$ & (b) $\intenergy(\sigma,\pi)$ 
\end{tabular}
}
\end{center}
\caption{\label{tablesix}
(a) the (symmetric) matrix illustrating $\destat(\sigma^{-1}\pi)$ for all $\sigma,\pi \in S_3$ where the rows and columns are indexed with the permutations from $S_3$ in lexicographic order.
(b) the (symmetric) matrix illustrating $\intenergy(\sigma,\pi):=1-\destat(\sigma^{-1}\pi)/2$ for all $\sigma,\pi \in S_3$ where the rows and columns are indexed with the permutations from $S_3$ in lexicographic order.
}
\end{table}
So if $\vec{\pi}$ is a configuration with $n_{\sigma}$ permaspins equal to $\sigma$, then the Hamiltonian $\H^{(mean)}(\vec{\pi})$ equals $\H^{(mean)} (n_{123},\ldots,n_{321})$ where
\begin{align*}
\lefteqn{\H^{(mean)} (n_{123},\ldots,n_{321}) }\\ 
\qquad =& - \dfrac{qJ}{n-1} 
\left( \sum_{\sigma \in S_3} \binom{n_{\sigma}}{2} \intenergy(\sigma,\sigma)  
 + \sum_{\sigma\prec \tau \in S_3} n_{\sigma}n_{\tau} \intenergy (\sigma,\tau) \right)
 -H \sum_{\sigma \in S_3} n_{\sigma} \intenergy(\sigma,\id)\\
=& - \dfrac{qJ}{n-1} 
\left( \sum_{\sigma \in S_3} \binom{n_{\sigma}}{2} 
 - (n_{123}n_{321}+n_{132}n_{231} +n_{213}n_{312})\right)
 -H\left( n_{123} - n_{321}\right).
\end{align*}
The binomial coefficients give rise to quadratic terms that can be neatly rewritten:
\begin{align*}
\lefteqn{\H^{(mean)} (n_{123},\ldots,n_{321}) }\\
=& - \tfrac{qJ}{2(n-1)}
\left( 
(n_{123}-n_{321})^2 + (n_{132}-n_{231})^2 + (n_{213}-n_{312})^2 - (n_{123}+\ldots+n_{321})
\right)\\ & \qquad  -H\left( n_{123} - n_{321}\right)\\
=& - \tfrac{qJ}{2(n-1)} 
\left( 
(n_{123}-n_{321})^2 + (n_{132}-n_{231})^2 + (n_{213}-n_{312})^2 - n
\right) 
 -H\left( n_{123} - n_{321}\right)\\
=& \tfrac{nqJ}{2(n-1)} 
+ \tfrac{H^2(n-1)}{2qJ}
- \tfrac{qJ}{2(n-1)} 
\left( 
\left(n_{123}-n_{321}+\tfrac{H(n-1)}{qJ}\right)^2 ~+~ (n_{132}-n_{231})^2 + (n_{213}-n_{312})^2  
\right).
\end{align*}
The partition function is 
\begin{align}
Z_n(\beta) =& \sum_{n_{123}+\ldots+n_{321}=n} \binom{n}{n_{123},n_{132},\ldots,n_{321}} 
	\exp \left(
	-\beta \H^{(mean)}(n_{123},\ldots,n_{321})
	\right). \label{mfpf}
\end{align}
Using the expression for the mean-field Hamiltonian above, we can now write the partition function as a reduced sum:

\begin{proposition}
\label{meanfieldsum}
The mean-field partition function is
\begin{align}
Z_n^{(mean)}(\beta) =& \exp\left({-\frac{\beta}{2}\left(\frac{nq J}{n-1} + \frac{ H^2(n-1)}{qJ}\right)}\right)\nonumber\\
& \quad
		\sum_{a+b+c =n}  \binom{n}{a,b,c} 
		G_a\left(\tfrac{H(n-1)}{qJ};e^{\frac{\beta qJ}{2(n-1)}}\right)
		G_b\left(0;e^{\frac{\beta qJ}{2(n-1)}}\right)
		G_c\left(0;e^{\frac{\beta qJ}{2(n-1)}}\right)\label{mfpf2}
\end{align}
where 
$$G_m(\ell;x) := \sum_{i+j=m} \binom{m}{i} x^{(i-j+\ell)^2} = \sum_{i} \binom{m}{i}  x^{(2i-m+\ell)^2}.$$
\end{proposition}
Given the quadratic term in the exponent of $x$, there is no closed formula for $G_m(\ell,x)$ that can help us achieve a closed form for $Z_n(\beta)$. 
This does not rule out the possibility of approximating $Z_n(\beta)$ for some special cases.
For the high-temperature case, which corresponds to $\beta$ being very small, let us observe that the exponential term in Equation~\ref{mfpf} is 
$$1-\beta  \H^{(mean)}(n_{123},\ldots,n_{321}) +O(\beta^2),$$ 
so the partition function for this case is given by a simple application of the multinomial theorem and we will have
$$Z_n^{(mean)}(\beta) \approx 6^n.$$
The same reasoning can be used to show for the high-temperature mean-field $k$-permaspin case we will have 
$$Z_n^{(mean)}(\beta) \approx (k!)^n.$$
These results can also be seen from the non-mean-field case as for $\beta$ small, the partition function can be written as a product of the partition function of the $n$ consistent subsystems (whose interaction energy is negligible), each of which has $k!$ different states.

Consider now the case of no external field, $H=0$, and when every vertex is a neighbour of every other vertex. This corresponds to $q=n-1$ which is essentially the case of the mean-field model on the complete graph.
Let $w=w(\beta)=e^{{\beta J}/{2}}$.
We have
\begin{align}\label{hiup}
Z_n^{(mean)}(\beta) =& w^{-n}
		\sum_{a+b+c =n}  \binom{n}{a,b,c} 
		g_a\left(w \right)
		g_b\left(w \right)
		g_c\left(w \right),
\end{align}
where 
$$g_m(x) :=  \sum_{i} \binom{m}{i}  x^{(2i-m)^2}.$$
Let $G(x,z) := \sum_{m \geq 0} g_m(x) z^m/m!.$
Multiplying both sides of Equation~\ref{hiup} by $z^n$ and reorganizing, the above equation becomes
\begin{align}
\dfrac{Z_n^{(mean)}(\beta)w^nz^n}{n!}  =& 
		\sum_{a+b+c =n}  \dfrac{1}{a!b!c!} 
		g_a\left(w \right) z^a
		g_b\left(w \right) z^b
		g_c\left(w \right) z^c.
\end{align}
By summing both sides over non-negative $n$, this allows us to give an exact expression for the exponential generating function of the mean-field partition function in terms of a power of a simple (but difficult to analyse) generating function involving the binomial coefficients.
$$
\sum_{n\geq 0} Z_n^{(mean)}(\beta) \dfrac{(wz)^n}{n!} 
= G^3(w,z).
$$
Note that $g_0(x)=1$, $g_1(x) = 2x$, $g_2(x) = 1+2x^4$, $g_3(x) = 2(x+x^9)$, $g_4(x) = 1+2(x^{4} +x^{16})$.
The largest term in $g_m(x)$ depends on whether $x \gtrless 1$.
If $x<1$ then the largest term in $g_m(x)$ will be attained when $2i-m$ is as small as possible. 
When $m$ is even this will be attained at $i=m/2$ and will be $\binom{m}{m/2}$ whereas for $m$ odd it will be $\binom{m}{(m+1)/2} x$.
If $x>1$ then the largest term in $g_m(x)$ will be attained when $2i-m$ is as large as possible. This occurs for $i=m$ and will be $x^{m^2}$.

%%%%%%%% %%%%%%%% %%%%%%%% %%%%%%%% %%%%%%%% %%%%%%%%
%%%%%%%% %%%%%%%% %%%%%%%% %%%%%%%% %%%%%%%% %%%%%%%%
%%%%%%%% %%%%%%%% %%%%%%%% %%%%%%%% %%%%%%%% %%%%%%%%

\section{Remarks}
\label{sectionseven}
This paper introduced a variant of the Ising model in which permutations and permutation statistics play a leading role.
This was motivated by a complexity measure for declarative systems.
While the 2-permaspin model corresponds exactly to the classical Ising model, we explored the next non-trivial case, the 1D 3-permaspin model.  
For the 1D 3-permaspin model in the absence of an external field, we were able to give exact expressions for the partition function and free energy based on the generating function of the $\cddes$ statistic, and were also able to give an expression for the 
free energy for the case when there are a large number of allowed spins. 

We saw that the calculations involved in the derivations of the partition function and free energy for the 1D 3-permaspin model in an external field were challenging and non-trivial.
Experimental evidence was used to provide several results for this case in relation to the free energy.
We also considered a mean-field version of the model and saw that the partition function the 1D 3-permaspin case could be neatly written as a certain trinomial sum. 
The no-external-field case admitted a particularly nice form in terms of an exponential generating function.

In our discussion in Section 2, we noted that the {\it{number of inversions}} permutation statistic could be another natural choice to replace $\cddes$ in the calculation of the free energy. 
The reason is that the $\inv$ statistic satisfies the symmetry property given in Equation~\ref{statsym} without modification. 
In order to consider the 1D 3-permaspin model using the $\inv$ statistic in place of $\cddes$, the transition matrix is:
$$
A^{(\mathrm{inv})}=\left(\begin{array}{rrrrrr}
1 & a b & a b & a^{2} b^{2} & a^{2} b^{2} & a^{3} b^{3} \\
a b & a^{2} & a^{2} b^{2} & a^{3} b^{3} & a^{3} b & a^{4} b^{2} \\
a b & a^{2} b^{2} & a^{2} & a^{3} b & a^{3} b^{3} & a^{4} b^{2} \\
a^{2} b^{2} & a^{3} b^{3} & a^{3} b & a^{4} & a^{4} b^{2} & a^{5} b \\
a^{2} b^{2} & a^{3} b & a^{3} b^{3} & a^{4} b^{2} & a^{4} & a^{5} b \\
a^{3} b^{3} & a^{4} b^{2} & a^{4} b^{2} & a^{5} b & a^{5} b & a^{6}
\end{array}\right).
$$
The characteristic polynomial for this matrix is
\begin{align*}
c_{A^{(\mathrm{inv})}}(\lambda) = &
\Big( (b^2 + b + 1)(b^2 - b + 1)a^{12}(b + 1)^4(b - 1)^4 \\
        & \quad - \lambda (a^4b^2 - a^2b^4 + a^4 + b^2 + 1)(a^2 + 1)a^6(b + 1)^2(b - 1)^2 \\
        & \quad - \lambda^2 (a^8 + 2a^6b^2 + 2a^4b^4 + a^6 + 3a^4b^2 + 2a^4 + 2a^2b^2 + a^2 + 1)a^2(b^2 - 1)\\
        & \quad +\lambda^3 (-(a^4 + a^2b^2 + 1)(a^2 + 1)) +\lambda^4\Big) \\
    &\Big(
    a^6(1+b)^3(1-b)^3 + \lambda (a^2 + 1)a^2(b + 1)(b - 1) +\lambda^2
    \Big).
\end{align*}
We note that the determination of the eigenvalues now relies on solving a quartic that does not factorize any further.

A main goal when dealing with permutation statistics is to derive their generating function since it encodes the distribution of the statistic for arbitrary length permutations. 
In this paper we have seen that information about such a generating function is useful for the zero-field case, but not so for the case of an external field.
Instead, determining the eigenvalues of a (transition) matrix whose exponents are encoded by the permutation statistic is the main goal. 
Might these eigenvalues have other uses, perhaps in some spectral theory of permutations that has yet to be formalized?

\end{document}